%% file: MR-necklace.tex
\documentclass{amsart}

\usepackage{graphicx,color,epsfig}

\usepackage{amssymb}
\theoremstyle{plain}

\newtheorem{thm}{Theorem}
\newtheorem*{thm*}{Theorem}
\newtheorem{lem}{Lemma}

\newtheorem{cor}{Corollary}
\newtheorem*{conj*}{Conjecture}
\newtheorem{prop}{Proposition}

\theoremstyle{definition}

\newtheorem{definition}{Definition}

\DeclareMathOperator{\verts}{vert}
\DeclareMathOperator{\supp}{supp}
\DeclareMathOperator{\tr}{tr}

\title[Splitting multidimensional necklaces]{Splitting multidimensional
necklaces}

\author{Mark de Longueville}
\author{Rade T.\ \v{Z}ivaljevi\'{c}}
\address{Freie Universit\"{a}t Berlin,
Fachbereich Mathematik, Arnimallee 3, 14195 Berlin, Germany}
\address{Mathematical Institute SANU, Knez Mihailova 35/1, p.f.\ 367, 11001 Belgrade, Serbia}
\date{October 2006}

\begin{document}
\begin{abstract}
The well-known ``splitting necklace theorem'' of Alon \cite{Alo87}
says that each necklace with $k\cdot a_i$ beads of color
$i=1,\ldots, n$ can be fairly divided between $k$ ``thieves'' by
at most $n(k-1)$ cuts. Alon deduced this result from the fact that
such a division is possible also in the case of a continuous
necklace $[0,1]$ where beads of given color are interpreted as
measurable sets $A_i\subset [0,1]$ (or more generally as
continuous measures $\mu_i$). We demonstrate that Alon's result is
a special case of a multidimensional, consensus division theorem
of $n$ continuous probability measures $\mu_1,\ldots ,\mu_n$ on a
$d$-cube $[0,1]^d$. The dissection is performed by $m_1+\ldots
+m_d=n(k-1)$ hyperplanes parallel to the sides of $[0,1]^d$
dividing the cube into $m_1\cdot\ldots\cdot m_d$ elementary
parallelepipeds where the integers $m_i$ are prescribed in
advance.
\end{abstract}
\maketitle

\section{Introduction}
\label{sec:intro}

The problem of {\em consensus division} arises when two or more
competitive or cooperative parties, each guided by their
individual objective functions, divide an object according to some
notion of fairness. There are many different mathematical
reformulations of this problem depending on what kind of divisions
are allowed, what kind of object is divided, whether the parties
involved are cooperative or not, etc. Early examples of problems
and results of this type are the ``ham sandwich theorem'' of
Steinhaus and Banach, the envy-free ``cake-division problem'' of
Steinhaus, the equipartition of  measurable sets by hyperplanes of
Gr\"{u}nbaum and Hadwiger, and more recently the ``splitting
necklace theorem'' of Alon,
\cite{Alo87,Gru60,{Hadw66},{Mat03},{Ziv04}}. A model example of a
fair-division theorem when two parties are involved is the
Hobby-Rice theorem.
\begin{thm}\label{thm:Hobby-Rice} {\rm (\cite{Ho-Ri65})}
Let $\mu_1,\mu_2,\ldots, \mu_n$ be a collection of continuous
probability measures on $[0,1]$. Then there exists a partition of
$[0,1]$ by $n$ cut points into $n+1$ intervals
$I_0,I_1,\ldots,I_n$ and the corresponding signs
$\epsilon_0,\epsilon_1,\ldots,\epsilon_n\in\{-1,+1\}$ such that
for each measure $\mu_i$,
$$
\sum_{j=0}^n \epsilon_j\cdot \mu_i(I_j) =0.
$$
\end{thm}
A well known consequence of this result is the ``necklace
theorem'', proved by Goldberg and West \cite{GoWe85}, which says
that every open necklace with $d$ kind of stones (an even number
of each kind) can be divided between two thieves using no more
than $d$ cuts.

A celebrated generalization of Theorem~\ref{thm:Hobby-Rice} is the
following ``splitting necklace theorem'' of Alon, which extends
the result of Goldberg and West to the case of $q$ ``thieves''. We
formulate the continuous version which includes
Theorem~\ref{thm:Hobby-Rice} as a special case and which can be
used to deduce the corresponding discrete version.

\begin{thm}\label{thm:Alon} {\rm (\cite{Alo87})} Let $\mu_1,\mu_2,\ldots,
\mu_n$ be a collection of $n$ continuous probability measures on
$[0,1]$. Let $k\geq 2$ and $N:=n(k-1)$. Then there exists a
partition of $[0,1]$ by $N$ cut points into $N+1$ intervals
$I_0,I_1,\ldots,I_N$ and a function $f : \{0,1,\ldots,
N\}\rightarrow \{1,\ldots, k\}$ such that for each $\mu_i$ and
each $j\in\{1,2,\ldots,k\}$,
$$
\sum_{f(p)=j} \mu_i(I_p) = 1/k .
$$
\end{thm}
Our main objective in this paper is to show that there exist
higher dimensional analogs (Theorems~\ref{thm:2-dim-necklace} and
\ref{thm:main}) of the splitting necklace theorem which include
Theorems~\ref{thm:Hobby-Rice} and \ref{thm:Alon} as special cases.
This may sound as a surprise in light of Theorem~5.2 from
\cite{AW86} claiming that, given $l\geq 0$, for every $d\geq 2$
there exist $2$-colorings of $[0,1]^d$ which do not admit
``bisections'' of size at most $l$. This ambiguity is immediately
resolved by the observation that the ``bisections'' allowed in
\cite{AW86} were of quite special nature ($d$-dimensional
checkerboards) while in our approach there are no restrictions on
the coloring (labelling) of elementary parallelepipeds.

An important step leading to the generalization of the ``splitting
necklace theorem'' was the recognition of the role of ``rainbow
complexes'' $\Omega(Q;S)$ where $Q$ is an arbitrary
$d$-dimensional, convex polytope and $S$ a finite set of
``colors'' used for labelling the vertices of $Q$. These complexes
turn out to be (topologically) shellable
(Theorem~\ref{thm:topo-shell}) and to have other interesting
properties reflecting the geometry and combinatorics of the base
polytope $Q$, Section~\ref{sec:concluding}.

\section{Two--dimensional necklaces and the configuration space $\Omega(m,n)$ }
\label{sec:2-dim-case}

As a preliminary step, before we address the general case of a
$d$-dimensional necklace ($d$-dimensional carpet) $I^d=[0,1]^d$,
with $n$ measures $\mu_1,\mu_2,\ldots,\mu_n$ on $I^d$, and $k$
parties (thieves) interested in a fair division, we focus our
attention on the case $d=k=2$. This case exhibits all the main
features of the general $d$-dimensional problem and provides a
motivation for the introduction of configuration spaces
$\Omega(m,n)$ and their generalizations.

A ``splitting'' of a square $I^2 = [0,1]\times [0,1]$ is a
partition of $I^2$ into smaller rectangles by lines parallel to
the sides of the square. Assuming that the square is positioned in
the coordinate system so that the diagonally opposite vertices are
$(0,0)$ and $(1,1)$, a $(m\times n)$-partition is determined by a
choice of $m$ points $0=x_0\leq x_1\leq x_2\ldots \leq x_m\leq
x_{m+1}=1$ on the $x$-axes and $n$ points $0=y_0\leq y_1\leq
y_2\ldots \leq y_n\leq y_{n+1}=1$ on the $y$-axes.

The associated splitting (partition) is the division of $I^2$ into
(possibly degenerate) rectangles $[x_i,x_{i+1}]\times
[y_j,y_{j+1}]$, where $i=0,\ldots, m$ and $j=0,\ldots , n$. Recall
an elementary fact that the space of all $m$-partitions $0=x_0\leq
x_1\leq x_2\ldots \leq x_m\leq x_{m+1}=1$ of the unit interval $I$
is naturally identified as the simplex $\Delta_m$ where $t_i :=
x_{i+1} - x_i$ are the associated barycentric coordinates.
Similarly, the barycentric coordinates associated to a
$y$-partition are $s_j = y_{j+1}-y_j$. It follows that the space
of all $(m\times n)$-partitions of the square $I^2$ is naturally
parameterized by points of the product $\Delta_m\times \Delta_n$.

\medskip
The basic cell $C_{(m,n)}=\Delta_m\times\Delta_n$ should play in
the case of $2$-dimensional partitions the role analogous to the
role the cell $C_m = \Delta_m$ plays in the case of
$1$-dimensional partitions. The next step is to introduce two
``thieves'' or players who want to divide among themselves
elementary rectangles $R_{(i,j)}=[x_i,x_{i+1}]\times
[y_j,y_{j+1}]$ arising from the subdivision. By construction, the
degenerate elementary rectangles, i.e.\ the rectangles such that
either $x_i=x_{i+1}$ or $y_j=y_{j+1}$ are allowed. However, it is
instructive to keep in mind that the ``thieves'' are primarily
interested in non-degenerated rectangles.

In the $1$-dimensional case, a division of intervals between two
thieves was described by a function $\omega : \widetilde{m}
\rightarrow \{+,-\}$ where $\widetilde{m}:=\{0,1,\ldots, m\}$ and
$\omega(i)= +$ (alternatively $\omega(i)=-$) means that the
interval $[x_i,x_{i+1}]$ was allocated to the first (respectively
second) player.

Similarly, in $2$-dimensions a function $\omega :
\widetilde{m}\times \widetilde{n}\rightarrow \{+,-\}$ completely
describes an allocation of elementary rectangles to the two
players.

As in the $1$-dimensional case, a natural {\em configuration space}
$\Omega(m,n)$ for the $2$-dimensional problem should take into account
all $(m\times n)$-partitions of $I^2$ together with all possible
allocation functions $\omega\in \{+,-\}^{\widetilde{m}\times
  \widetilde{n}}$. In other words a typical element in $\Omega(m,n)$
is a triple $(t,s;\omega)\in C_{(m,n)}\times
\{+,-\}^{\widetilde{m}\times \widetilde{n}}$.  Collecting together
all triples $(t,s;\omega)$ corresponding to a fixed $\omega\in
\{+,-\}^{\widetilde{m}\times \widetilde{n}}$ we observe that
$\Omega(m,n)$ ought to be the union of cells $C_{(m,n)}^\omega:=
\{(t,s;\omega)\mid (t,s)\in\Delta_m\times\Delta_n\}$. Two cells
$C_{(m,n)}^\omega$ and $C_{(m,n)}^\nu$ can have a point in common.
This happens precisely if whenever $\omega(i,j)\neq \nu(i,j)$, the
corresponding rectangle $R_{(i,j)}=[x_i,x_{i+1}]\times
[y_j,y_{j+1}]$ is degenerate. This leads us to the definition of
the following space

$$
\Omega(m,n)=\coprod_{\omega\in \{+,-\}^{\widetilde{m}\times \widetilde{n}}}
C_{(m,n)}^\omega /\approx
$$
where $(t,s;\omega)\approx (t',s';\omega')$ if and only if $t=t'
\mbox{ {\rm and} } s=s'$ and
$$(\omega(i,j)\neq \omega'(i,j)) \Rightarrow (t_i=t'_i=0 \mbox{ {\rm or} }
s_j=s'_j=0).$$

\begin{figure}[hbt]
\centering
\includegraphics[scale=0.35]{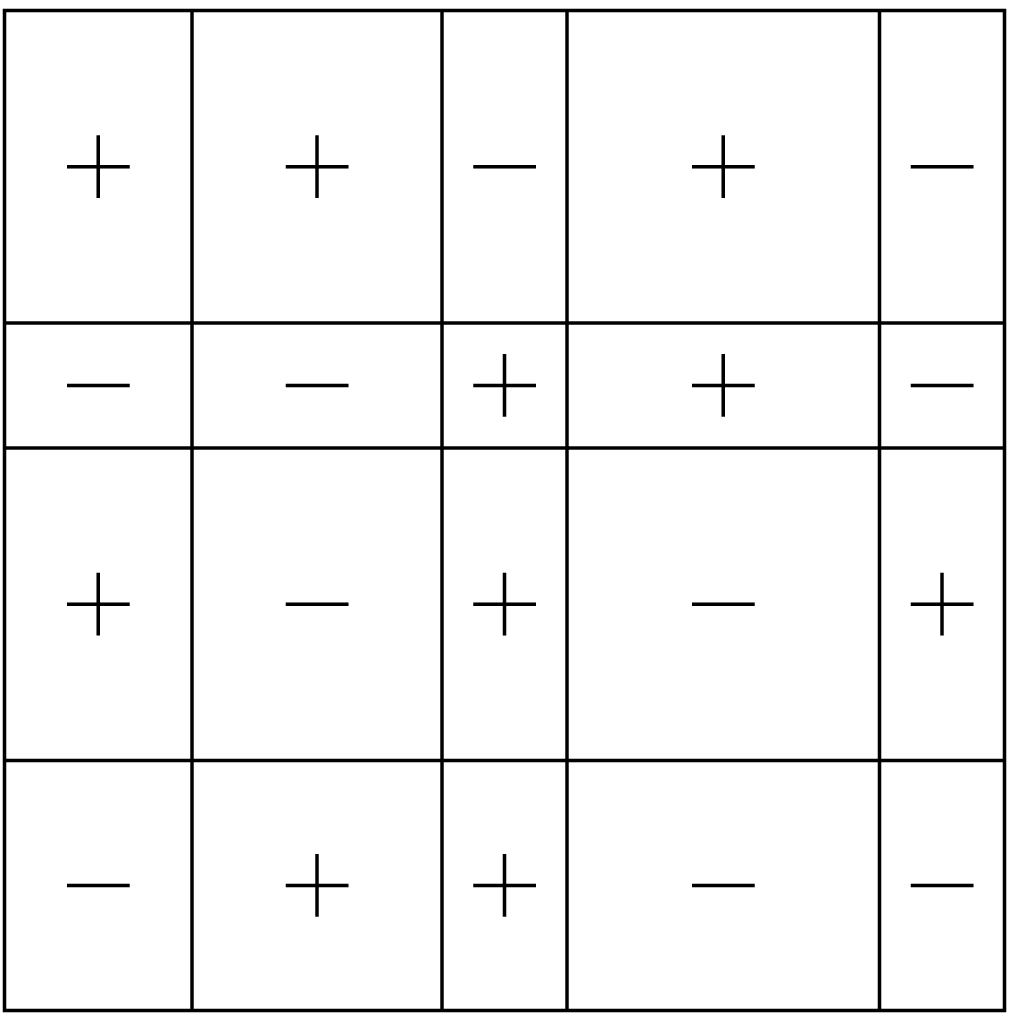}
\caption{An element of $\Omega(m,n)$.} \label{fig:chessboard}
\end{figure}

Here is a convenient way to ``visualize'' the configuration space
$\Omega(m,n)$. An element $x=(t,s;\omega)\in\Omega(m,n)$ is
visualized as a $(m+1)\times(n+1)$-``chessboard'' where the pair
$(t,s)\in\Delta_m\times\Delta_n$ determines the size and the shape
of each of the elementary parallelepipeds while the coloring
(labelling) is described by the function $\omega$
(Figure~\ref{fig:chessboard}).

\begin{figure}[htbp]
  \centering
  \input{omegacell.pstex_t}
  \caption{A part of $\Omega(1,1)$ and a labelled cell of $\Omega(2,1)$}
  \label{fig:omegacell}
\end{figure}
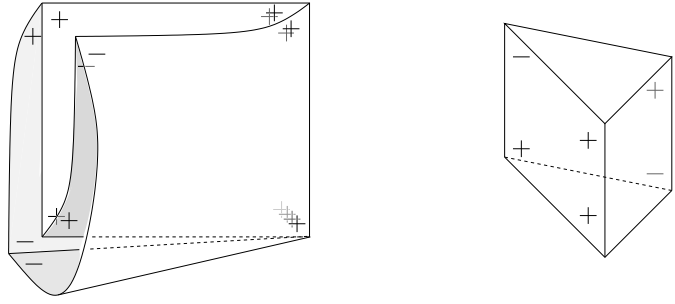

Each cell $C_{(m,n)}^\omega$ is visualized as the polytope
$C_{(m,n)}:=\Delta_m\times\Delta_n$ with vertices colored
(labelled) by $+$ or $-$, according to the prescription given by
$\omega$, while the total configuration space $\Omega(m,n)$ is the
union of cells $C_{(m,n)}^\omega$ (Figure~\ref{fig:omegacell}).
Note that the elementary parallelepipeds from
Figure~\ref{fig:chessboard} are in one-to-one correspondence with
the vertices of $C_{(m,n)}$ so one can read off the labelling
function $\omega$ both from the coloring of the elementary
parallelepipeds and the coloring of the vertices of the cell
$C_{(m,n)}$.

The proof of the two--dimensional analogue of Alon's theorem
relies on the following important property of the configuration
space $\Omega(m,n)$.
\begin{thm}
  The configuration space $\Omega(m,n)$ is a
$(m+n)$-dimensional, $(m+n-1)$-connected, free $\mathbb{Z}_2$-complex.
\label{thm:connectivity-2-case}
\end{thm}

In subsequent sections we will obtain a stronger and much more
general result. However, here we present an outline of a direct
proof  of this theorem which provides additional insight into the
structure of complexes $\Omega(m,n)$.

\begin{proof}[Sketch of proof]
We proceed  by induction on $\nu=m+n$. The complexes $\Omega(m,0)$
and $\Omega(0,n)$ are isomorphic to $[2]^{\ast(m+1)}\cong S^m$ and
$[2]^{\ast(n+1)}\cong S^n$ respectively. The complex
$[2]^{\ast(m+1)}=[2]\ast\ldots\ast[2]$ is naturally isomorphic to the
boundary $\partial\Diamond_m$ of the crosspolytope $\Diamond_m :={\rm
  conv}\{e_i,-e_i\}_{i=1}^m\subset \mathbb{R}^m$. This holds also in
the case $m=n=0$, i.e.\ $\Omega(0,0)\cong S^0\cong [2]$ is the
boundary of $\partial\Diamond_1\cong [-1,+1]$.
\par
Surprisingly enough, the complex $\Omega=\Omega(m,n)$ exhibits
formal structure similar to the complex $\partial\Diamond_m$ in
the sense that it can be associated a ``north and south pole'' and
the ``upper and lower hemisphere'' $\Omega^+$ and $\Omega^-$ with
all the usual consequences including the associated Mayer-Vietoris
exact sequence of the triple $(\Omega;\Omega^+,\Omega^-)$.

In order to define a ``north and south pole'' in $\Omega(m,n)$,
let us start with a maximally degenerated partition
$(\bar{t},\bar{s})\in C_{(m,n)}$ where $\bar{t}_0=\bar{s}_0=1$ and
$\bar{t}_i=\bar{s}_j=0$ for $i\geq 1$ and $j\geq 1$. In this
partition there is only one non-degenerated elementary rectangle,
consequently there are only two associated elements $c_+:=
(\bar{t},\bar{s};+)$ and $c_-:= (\bar{t},\bar{s};-)$ in
$\Omega(m,n)$.

The ``upper hemisphere'' $\Omega^+(m,n)$ is the set of all points
$x=(t,s;\omega)\in \Omega(m,n)$ which are visible from $c_+$,
i.e.\ such that both $x$ and $c_+$ belong to the same maximal cell
in $\Omega(m,n)$. In other words $x\in\Omega^+(m,n)$ has a
representative $x=(t,s;\omega)$ such that $\omega(0,0)=+$. The
``lower hemisphere'' $\Omega^-(m,n)$ is defined similarly as the
set of all points $x=(t,s;\omega)$ which allow a representation
such that $\omega(0,0)=-$.

Both $\Omega^+(m,n)$ and $\Omega^-(m,n)$ are contractible. Indeed,
both spaces are star-shaped,  with centers $c_+$ and $c_-$
respectively, and a contraction is defined by the linear homotopy.

Let us focus on the structure of the ``equatorial set''
$E(m,n):=\Omega^+(m,n)\cap\Omega^-(m,n)$. By definition $x=
(t,s;\omega)\in E(m,n)$ if either $t_0=0$ or $s_0=0$. From here it
follows that $E(m,n)=A\cup B$ where $A\cong \Omega(m-1,n)$ and $B\cong
\Omega(m,n-1)$. Since $x\in A\cap B$ if and only if $t_0=s_0=0$, we
observe that $A\cap B\cong\Omega(m-1,n-1)$.  A twofold application of
the Mayer-Vietoris sequence to the triads
$(\Omega(m,n);\Omega^+(m,n),\Omega^-(m,n))$ and $(E(m,n);A,B)$
together with a Seiffert--van Kampen argument for determining the
fundamental group of $\Omega(m,n)$ yields the desired connectivity.
\end{proof}

Theorem~\ref{thm:connectivity-2-case}, following the usual
Configuration space/Test map scheme \cite{Ziv04}, is the basis for
the following version of the two--dimensional splitting necklace
theorem.

\begin{thm}[Two--dimensional necklace  for two thieves]
\label{thm:2-dim-necklace} Let $\mu_1,\ldots,\mu_n$ be a
collection of $n$ continuous probability measures on the unit
square $I^2=[0,1]^2$.  Then for any choice of $m_1,m_2\geq 0$ of
integers such that $m_1+m_2=n$, there exist $m_1$ vertical and
$m_2$ horizontal cuts of the square, and a coloring of the
elementary rectangles obtained this way by two colors ``$+$'' and
``$-$'' such that $\mu_i(A_+)=\mu_i(A_-)=\frac{1}{2}$ for all $i$
where $A_+$ (respectively $A_-$) is the union of all elementary
parallelepipeds colored by ``$+$'' (respectively ``$-$'').
\end{thm}

We omit the proof of Theorem~\ref{thm:2-dim-necklace} since it
will be subsumed by a more general argument used in the proof of
Theorem~\ref{thm:main} and instead turn our attention to the
general case of a necklace in $d$ dimensions for an arbitrary
number of thieves.

%The
%Mayer-Vietoris exact sequence of the triad
%$(\Omega(m,n);\Omega^+(m,n),\Omega^-(m,n))$ splits into sequences
% \begin{equation}\label{eqn:M-V-1}
% 0 \longrightarrow
%H_{k+1}(\Omega(m,n))\longrightarrow H_k(E(m,n))\longrightarrow 0
% \end{equation}
%which implies the isomorphism $H_{k+1}(\Omega(m,n))\cong
%H_k(E(m,n))$ for each $k\geq 0$.
%The Mayer-Vietoris exact sequence
%of the triad $(E(m,n);A,B)$,
%\begin{equation}\label{eqn:M-V-1}
%H_{k}(\Omega(m-1,n))\oplus H_{k}(\Omega(m,n-1))\longrightarrow
%H_k(E(m,n))\longrightarrow H_{k-1}(\Omega(m-1,n-1))
% \end{equation}
%together with the induction hypothesis, allows us to deduce that
%$\widetilde{H}_j(E(m,n))=0$ for $j\leq m+n-2$. This in turn lead
%us to the conclusion that $\widetilde{H}_j(\Omega(m,n))=0$ for
%$j\leq m+n-1$.

%\medskip
%In order to complete the proof it is sufficient to show that
%$\pi_1(\Omega(m,n))=0$ for

\section{The complex $\Omega(Q;G)$ of $G$-labelled polytopes}
\label{sec:G-labelled polytopes}

The $2$-dimensional splitting necklace theorem presented in
Section~\ref{sec:2-dim-case}, especially the construction of the
configuration space $\Omega(m,n)$ with favorable properties
(Theorem~\ref{thm:connectivity-2-case}), reveal that higher
dimensional analogs and extensions should be within reach by
similar methods. Apparently the most natural generalization that
comes to mind is the splitting of a $d$-dimensional cube by
hyperplanes parallel to its sides. Moreover, in order to extend
the $1$-dimensional ``splitting necklace theorem'', we should
replace $\{+,-\}$ by an arbitrary set $G$ of labels (colors)
corresponding to different ``thieves''. The letter $G$ should
indicate that the labels are often elements of a given finite
group, e.g.\ $G\cong \mathbb{Z}_2\cong\{+,-\}$ in the case of two
thieves.

\par
An extension and a multidimensional analogue of the configuration
space $\Omega(m,n)=\Omega(m,n;\pm)$ is the space
$\Omega(\mathfrak{m};G)$ = $\Omega(m_1,m_2,\ldots, m_d;G)$ defined
as follows. A typical element in $\Omega(\mathfrak{m},G)$ is a
pair $(\mathfrak{t};\omega)\in Q_{\mathfrak{m}}\times
G^{\widetilde{m}_1\times\ldots\times \widetilde{m}_d}$ where
$Q_{\mathfrak{m}}:=\Delta_{m_1}\times\ldots\times \Delta_{m_d}$ is
the space of all $\mathfrak{m}$-partitions of the hypercube $I^d$.
More precisely each of the coordinates $t_i$ of
$\mathfrak{t}=(t_1,\ldots, t_d)$ is a partition $0=x_0^i\leq
x_1^i\leq\ldots \leq x_{m_i}^i\leq x_{m_i+1}^i\leq 1$ of the
interval $[0,1]$ so an elementary (possibly degenerate)
$d$-parallelepiped associated to $\mathfrak{t}$, indexed by
$\mathfrak{j}=(j_1,\ldots, j_d)\in
\widetilde{m}_1\times\ldots\times \widetilde{m}_d$, is
$$
R_{\mathfrak{j}}(\mathfrak{t})=[x_{j_1}^1,x_{j_1+1}^1]\times\ldots
\times [x_{j_d}^d,x_{j_d+1}^d].
$$
For each labelling function $\omega :
\widetilde{m}_1\times\ldots\times \widetilde{m}_d \rightarrow G$
there is an associated cell $C^\omega_{\mathfrak{m}}\cong
\Delta_{m_1}\times\ldots\times \Delta_{m_d}$. It is convenient to
visualize the cell $C^\omega_{\mathfrak{m}}$ as the polytope
$Q=\Delta_{m_1}\times\ldots\times \Delta_{m_d}$ with all vertices
labelled by elements from $G$. This leads us to the following
definition.

\begin{definition}\label{def:omega-G-1} The configuration space
$\Omega(\mathfrak{m};G)$ is defined as the quotient space:
$$
\coprod_{\omega\in {G}^{\widetilde{m}_1\times\ldots\times
\widetilde{m}_d}} C_{\mathfrak{m}}^\omega /\approx
$$
where $(\mathfrak{t};\omega)\approx (\mathfrak{s};\nu)$ if and
only if $\mathfrak t = \mathfrak s$ and
$$\omega(\mathfrak{j})\neq \nu(\mathfrak{j})
\Rightarrow R_{\mathfrak{j}}(\mathfrak{t})=
R_{\mathfrak{j}}(\mathfrak{s}) \mbox{ {\rm is a degenerated
$d$-parallelepiped.} }
$$
\end{definition}

%\begin{definition}
%  Let $G$ be a group, and $m_1,\ldots,m_d\geq 0$. Define
%  $\Omega(m_1,\ldots,m_d;G)$ to be the following configuration space.
%\begin{align*}
%  \Omega(m_1,\ldots,m_d;G)=\left\{(s_1,\ldots,s_d,\omega): s_i\in
%    \Delta_{m_i},\omega\in G^{\{0,\ldots,m_1\}\times\cdots\times\{0,\ldots,m_d\}}\right\}/\sim
%\end{align*}
%We think of $\omega$ as being a labelling
%$\omega:\verts(\prod_{i=1}^d\Delta_{m_i})\rightarrow G$. The gluing
%$\sim$ is defined as follows. Let
%\begin{align*}
%  \supp(s_1,\ldots,s_d)&=\left\{(j_1,\ldots,j_d):(s_i)_{j_i}\not=0\text{ for all }i\right\}\\
%&\subseteq \{0,\ldots,m_1\}\times\cdots\times\{0,\ldots,m_d\}.
%\end{align*}
%and define
%\begin{align*}
%  (s_1,\ldots,s_d,\omega)\sim (t_1,\ldots,t_d,\tau)
%\end{align*}
%if and only if
%\begin{itemize}
%\item $s_i=t_i$ for all $i$, and
%\item if $(j_1,\ldots,j_d)\in\supp (s_1,\ldots,s_d)$, then $\omega(j_1,\ldots,j_d)=\tau(j_1,\ldots,j_d)$.
%\end{itemize}
%\end{definition}

%The construction of the space $\Omega(m_1,\ldots,m_d;G)$ can be
%generalized to arbitrary convex polytopes.

A natural extension of the configuration space
$\Omega(\mathfrak{m};G)$ is the cell complex $\Omega(Q;G)$ where
$Q$ is an arbitrary convex polytope $Q\subset \mathbb{R}^d$. Given
a function $\omega : \verts(Q)\rightarrow G$, the associated cell
$Q^\omega$ is described as the polytope with each vertex $v$
decorated (labelled) by the corresponding element $\omega(v)$. In
particular $Q^\omega=C_{\mathfrak{m}}^\omega$ if
$Q=\Delta_{m_1}\times\ldots \times\Delta_{m_d}$. Given
$\mathfrak{t}\in Q$, the associated element in $Q^\omega$ will be
denoted by $(\mathfrak{t},\omega)$. The cell $Q^\omega$ is
sometimes referred to as a {\em vertex-colored polytope} and
$\Omega(Q;G)$ is the associated {\em rainbow complex}.
\begin{definition}\label{def:omega-G-2} The configuration space
$\Omega(Q;G)$ is defined as the quotient space:
$$
\coprod_{\omega\in G^{\verts(G)}} Q^\omega /\approx
$$
where $(\mathfrak{t},\omega)\approx (\mathfrak{s},\nu)$ if and
only if $\mathfrak t = \mathfrak s$ and if $F\subset Q$ is the
minimal face such that $\mathfrak t\in F$, then
$\omega|_{\verts(F)}=\nu|_{\verts(F)}$.
\end{definition}

%\begin{definition}
%Let $P$ be a convex polytope with vertex set $V$ and $S$ be any set.
%For $x\in P$ let $F$ be the smallest face of $P$ containing $x$ and
%define $\supp(x)=V\cap F$, i.e., the vertices contained in $F$.
%Define $\Omega(P;S)$ by
%\begin{align*}
%  \Omega(P;S)=\{(x,\omega):x\in P, \omega:V\rightarrow S\}/\sim,
%\end{align*}
%where $(x,\omega)\sim(y,\tau)$ if and only if $x=y$ and $\omega|_{\supp(x)}=\tau|_{\supp(x)}$.
%\end{definition}

%$\Omega(\frak m;G)=\Omega(m_1,\ldots,m_d;G)$ is obtained in this way by setting
%$Q=\prod_{i=1}^d\Delta_{m_i}$.

\section{Shellability of $\Omega(Q;G)$}

One of the key ingredients in the proof of the higher dimensional
splitting necklace theorem is the proof that the complex
$\Omega(\mathfrak{m},G)$ is $(\vert\mathfrak{m}\vert-1)$-connected
where $\vert\mathfrak{m}\vert:= m_1+\ldots +m_d$. This could be
proved along the lines of the proof of
Theorem~\ref{thm:connectivity-2-case}. In this section we offer a
different proof of a more general fact that $\Omega(Q;G)$ is
always a $d$-dimensional, $(d-1)$-connected regular cell complex.

\subsection{Topological shellability}

A convenient way to prove that a (regular, polyhedral, simplicial)
$d$-dimensional cell complex is $(d-1)$-connected is to show that
it is {\em shellable} \cite{Bjo91} \cite{Zie95}. There are many
different concepts of shellability. Here, as a variation on a
theme, we introduce a form of shellability which will be referred
to as {\em topological shelling}.

\begin{definition}\label{def:topo-shell}
Suppose that $K$ is a finite, regular cell complex. A total
ordering $C_1,C_2,\ldots, C_k$ of its maximal cells is a {\em
topological shelling} of $K$ if
 $$
 {\rm dim}(C_1)\geq {\rm dim}(C_2)\geq \ldots \geq {\rm dim}(C_k)
 $$
and for each $j>1$ either $(a)$ or $(b)$ is satisfied where
\begin{enumerate}
\item[{\rm ($a$)}] $(\cup_{i<j} C_i)\cap C_j$ is a (non-empty)
contractible subset of $\partial(C_j)$,

\item[{\rm ($b$)}] $(\cup_{i<j} C_i)\cap C_j=\partial(C_j)$ where
$\partial(C_j)\cong S^{{\rm dim}(C_j)-1}$.
\end{enumerate}
\end{definition}
The following result is easily established by induction on the
number of maximal cells in $K$.
\begin{prop}\label{prop:topo-shell}
Suppose that $K$ is a finite cell complex which admits a
topological shelling $C_1,C_2,\ldots,C_k$. Let $n_i:={\rm
dim}(C_i)$. Then $K$ is homotopic to the wedge $\bigvee_{j\in
S}~S^{n_j}$ where $S:=\{j\mid (\cup_{i<j} C_i)\cap C_j=
\partial(C_j)\}.$
\end{prop}

\begin{proof}
Let $K_{\leqslant j} $ and $K_{< j}$
be the subcomplexes of $K$ defined by  $K_{\leqslant
j}:=\cup_{i\leqslant j}~C_i$ and $K_{< j}:=\cup_{i< j}~C_i$.
Suppose that by induction hypothesis the statement is true for all
$j< j_0$. If $K_{< j_0}\cap C_{j_0}$ is a contractible subset of
$\partial(C_{j_0})$ then $K_{<j_0}$ and $K_{\leqslant {j_0}}$ have
the same homotopy type, consequently $K_{\leqslant j_0}$ is also a
wedge of spheres.

Suppose $K_{< j_0}\cap C_{j_0}=\partial(C_{j_0})$. By the
induction hypothesis $K_{< j_0}$ is a wedge of spheres, $K_{<
j_0}\simeq \bigvee_{s=1}^t~S^{p_s}$ where
$$p:={\rm
min}\{p_s\}_{s=1}^t\geq {\rm dim}(C_{j_0})> {\rm
dim}(\partial(C_{j_0})).$$ It follows that $\partial(C_{j_0})$ is
contractible in $K_{< j_0}$, hence
$$
K_{j_0}\simeq S^{{\rm dim}(C_{j_0})}\vee \bigvee_{s=1}^t~S^{p_s} .
$$
\end{proof}

\begin{cor}\label{cor:connectedness}
A cell complex $K$ admitting a topological shelling is
$(n-1)$-connected, provided it is pure $n$-dimensional, i.e.\ if
all its maximal cells have the same dimension $n$.\qed
\end{cor}
\subsection{Topological shellability of $\Omega(Q;[k])$}

\begin{thm}\label{thm:topo-shell}
The complex $\Omega(Q;G)$ admits a topological shelling for each
convex $d$-polytope $Q\subset \mathbb{R}^d$ and each finite set $G$ of
labels (colors).
\end{thm}

\begin{proof}
Suppose that $k :=\vert G\vert$ is the cardinality of
the set $G$. If the polytope $Q = \Delta = \Delta_\nu$ is a
$\nu$-dimensional simplex then
$$\Omega(Q;G)\cong \Omega(\Delta_\nu;[k])=[k]\ast\ldots\ast [k]=[k]^{\ast
(\nu+1)}$$ is a simplicial complex which is well known to be
(lexicographically) shellable. Indeed, each $\nu$-dimensional
simplex in $[k]=[k]^{\ast (\nu+1)}$ is obtained from the simplex
$\Delta$ by coloring its vertices with colors from $[k]$. In other
words each of these simplexes is a vertex-colored polytope
$\Delta^f$ where $f: \{0,1,\ldots,\nu\}\rightarrow [k]$ is the
associated coloring function. Given functions $f,g\in
[k]^{\widetilde{\nu}}$, the lexicographical ordering defined by
\begin{equation}\label{eqn:lex}
f \prec g  \Leftrightarrow f(i)< g(i) \mbox{ {\rm where} } i:={\rm
min}\{j\mid f(j)\neq g(j)\}
\end{equation}
induces a shelling $\{\Delta^f\}_{f\in [k]^{\widetilde\nu}}$ of
the rainbow complex $\Omega(\Delta;[k])$. Indeed, for each $g\in
[k]^{\widetilde\nu}$ the complex $(\cup_{f\prec
g}~\Delta^f)\cap\Delta^g$ is easily shown to be a union of facets
of the simplex $\Delta^g$.

A convex polytope $Q\subset \mathbb{R}^d$ with vertices
$\verts(Q)=\{v_0,v_1,\ldots,v_\nu\}$ is the image ${\rm Im}(h)$ of
an associated affine map $h : \Delta_\nu\rightarrow Q, \, j\mapsto
v_j$. Given a function $f : \widetilde{m}\rightarrow [k]$, there
is an induced map $h^f : \Delta^f\rightarrow Q^f$ of
vertex-colored polytopes. This map however does not extend to a
cellular map of complexes $\Omega(\Delta;[k])$ and $\Omega(Q;[k])$
since the intersection $\Delta^f\cap\Delta^g$ is not necessarily
mapped to $Q^f\cap Q^g$. Nevertheless, the following claim shows
that both complexes admit formally the same shelling order.

\medskip\noindent
{\bf Claim:} The ordering $\{Q^f\}_{f\in [k]^{\widetilde{m}}}$
arising from the lexicographical ordering of functions
(\ref{eqn:lex}) is a topological shelling of the complex
$\Omega(Q;[k])$.

\medskip\noindent
{\bf Proof of the Claim:} Given a function $g\in
[k]^{\widetilde{m}}$, let $L_g$ be the complex
$$L_g:=\Omega_{\prec g}\cap Q^g=(\cup_{f\prec g} Q^f)\cap Q^g =
\cup_{f\prec g} (Q^f\cap Q^g).$$ According to
Definition~\ref{def:topo-shell} we have to demonstrate that $L_g$
is either contractible or $L_g=\partial(Q^g)\cong S^{d-1}$. The
intersection $Q^f\cap Q^g$ is the union of all vertex-colored
polytopes $F^h$ where $F$ is a face of $Q$ and $h :
\verts(F)\rightarrow [k]$ agrees with both $f$ and $g$ on
$\verts(F)$, i.e., $h=f\vert_{\verts(F)} = g\vert_{\verts(F)}$. In
the special case when $f(j)=g(j)$ for all but one element
$j_0\in\widetilde{\nu}$, i.e., if $\{j\in \widetilde\nu\mid
f(j)=g(j)\}=\widetilde\nu\setminus \{j_0\}$, we observe that
$Q^f\cap Q^g$ is essentially the ``anti-star'' $a$-$Star(v_{j_0})$
of the vertex $v_{j_0}$ in $\partial(Q^g)\cong\partial(Q)$, i.e.,
the union of all facets in $Q$ that do not contain the vertex
$v_{j_0}$. The anti-star corresponds to the facet
$\Delta^f\cap \Delta^g$ of $\Delta^f$, resp. $\Delta^g$, in the original shelling.

Given a face $F$ of
$Q$, let $open$-$Star(F)$ be the union of all relative interiors
of all proper faces of $Q$ which contain $F$ as a face
$$
open\mbox{-}Star(F) = \bigcup_{F\subseteq G\neq Q}
rel\mbox{-}int(G).
$$
One easily checks that
$$
a\mbox{-}Star(v) = \partial(Q)\setminus open\mbox{-}Star(v).
$$
In light of the fact that ``$\prec$'' is a shelling order of the
simplicial complex $[k]^{\ast (\widetilde\nu)}$,  i.e.,  $(\cup_{f\prec
g}~\Delta^f)\cap\Delta^g$ is a union of facets,
we observe that
there exists a non-empty set $S\subset\widetilde\nu$ such that
$$
L_g=\bigcup_{j\in S}~a\mbox{-}Star(v_j)=\partial(Q)\setminus
\bigcap_{j\in S} open\mbox{-}Star(v_j).
$$
Finally,
$$
\bigcap_{j\in S} open\mbox{-}Star(v_j) = open\mbox{-}Star(F)
$$
where $F:=\supp\{\{v_j\}\}_{j\in S}$ is the minimal face of
$Q$ containing all vertices $v_j$. It follows that
$$
L_g = \partial(Q)\setminus open\mbox{-}Star(F)
$$
which completes the proof since if the open star of $F$ is
non-empty, its complement is homeomorphic to a $(d-1)$-dimensional
cell.
\end{proof}

\section{The necklace theorem in arbitrary dimension}
%\begin{thm}
%  Let $m_1,\ldots,m_d\geq0$ and $m=m_1+\cdots+m_d$, then
%  $\Omega(m_1,\ldots,m_d;G)$ is an $E_mG$ space.
%\end{thm}

%\begin{thm}[\"Ozaydin, Sarkaria, Volovikov]
%  Let $p$ be prime, $n,r\geq 1$, and $G=\mathbb (Z_p)^r$ the $r$-th
%  power of the cyclic group, i.e., a
%  group of order $p^r$ and $N=n(p^r-1)$. Let $\Omega$ be any
%  $(N-1)$-connected free $G$-space. If $\mathbb E$ is linear
%  representation of $G$ of dimension $N$ with $\mathbb E^G=\{0\}$,
%  then every continuous $G$-map $f:\Omega\rightarrow \mathbb E$ has a zero.
%\end{thm}

%As we will need this result only in the case $G=\mathbb Z_p$ we give a
%proof of this case for the sake of its beauty and self containment of
%the article

The following result of Borsuk-Ulam type is a key tool for many
applications of  equivariant topological methods in combinatorics
and discrete geometry, \cite{Mat03,Ziv04}.
\begin{thm}[B\'ar\'any, Schlosman, Sz\"ucs \cite{bss}; Dold \cite{dold}]
\label{thm:BSS-D} Let $G=\mathbb{Z}_p$ be the cyclic group of
prime order $p$. Suppose that $\Omega$ is a finite,
$(N-1)$-connected, free $G$-cell complex where $N=n(p-1)$ for some
integer $n\geq 1$. Assume that $\mathbb E$ is a real, linear
$G$-representation of dimension $N$, having no trivial
subrepresentations, i.e.\ such that $\mathbb E^G=\{0\}$. Then
every continuous $G$-equivariant map $f:\Omega\rightarrow \mathbb
E$ has a zero.
\end{thm}
Although the proofs of this result are nowadays readily available
\cite{Mat03}, for the reader's convenience and self containment of
the paper we outline a short proof of this fact.
\begin{proof}
Assume that there is a map $f:\Omega\rightarrow \mathbb E$ without
a zero. This yields a $G$-equivariant map $\bar
f:\Omega\rightarrow S(\mathbb E)$ to the $(N-1)$-sphere
$S(\mathbb{E})$ in $\mathbb E$. As $p$ is prime and $0$ is the
only element in $\mathbb{E}$ fixed by all elements in $G$, it
follows that the induced action on $S(\mathbb E)$ is free. Hence
by the $(N-1)$-connectedness of $\Omega$ there exists a
$G$-equivariant map $g:S(\mathbb E)\rightarrow \Omega$. Now
consider the map $(g\circ \bar f)_\#:C_*(\Omega)\rightarrow
C_*(\Omega)$ for the cellular chain complex with respect to a
finite $G$-invariant cell structure. As every orbit of a cell
consists of $p$ elements, the Lefshetz trace $\Lambda(g\circ \bar
f)=\sum (-1)^i\tr (g\circ \bar f)_\#$ will be divisible by $p$. If
we compute the Lefshetz trace now on the homology level, we obtain
$\Lambda(g\circ \bar f)=\sum (-1)^i\tr (g\circ \bar f)_\ast=1$ as
the map factors through the homology of an $(N-1)$-sphere. A
contradiction!
\end{proof}

Given a set $\mathcal X$ of hyperplanes in $I^d$, let $C(\mathcal
X)$ be the set of cells (connected components) of $I^d\setminus
\bigcup \mathcal X$. For any coloring (labelling) map
$\omega:C(\mathcal X)\rightarrow [k]$, let $A_i:=\bigcup\{c\in
C(\mathcal X)\mid \omega(c)=i\}=\bigcup\omega^{-1}(i)$ be the
union of all cells colored by the same color $i$.
\begin{thm}[Higher dimensional necklace theorem]
\label{thm:main}
%d dimension of the space
%n number of measures
%k number of pieces
Assume $n,d\geq 1$ and $k\geq 2$, and let $\mu_1,\ldots,\mu_n$ be
a collection of $n$ continuous probability measures on the
$d$-dimensional cube $I^d\subset \mathbb{R}^d$. For any selection
of non-negative integers $m_1,\ldots,m_d$ such that
$m_1+\cdots+m_d=n(k-1)$ there exists a fair division with $m_i$
hyperplane cuts parallel to $i$-th coordinate hyperplane. In other
words there exists a set ${\mathcal X}=\bigcup_{i=1}^d{\mathcal
X}_i$ of $n(k-1)$ hyperplanes such that $\vert{\mathcal X}_i\vert
=m_i$, each $H\in {\mathcal X}_i$ is perpendicular to $e_i$, and
for some coloring function $\omega:C(\mathcal X)\rightarrow [k]$
$$
\mu_i(A_j)=\frac{1}{k} \mbox{ {\rm for all} }
i=1,\ldots,n \mbox{ {\rm and} } j=1,\ldots,k
$$
where $A_i:=\bigcup\omega^{-1}(i)$ are the unions of all cells colored by
the same color.
\end{thm}
We reduce the proof of the theorem to the
case $k=p$, $p$ prime.
\begin{lem}
  If the previous theorem holds for parameters $k_1,k_2\geq 2$ (in
  place of $k$) then it also holds for $k=k_1k_2$.
\end{lem}
\begin{proof}
  Before commencing the proof, the reader is referred to Figure~\ref{fig:induction}
  for a rough idea how the reduction claimed in the lemma is achieved.
  This is an example with $n(k-1)=6$ cutting hyperplanes where
  $n=2$, $d=2$, $k=k_1\cdot k_2=2\cdot2=4$, $m_1=2$, and $m_2=4$. The
  densities of the two measures $\mu_1$ and $\mu_2$ are indicated by the light and dark
  grey regions. The cube will be divided in the first step into $k_1=2$
  pieces.  Then the two pieces will be treated separately.
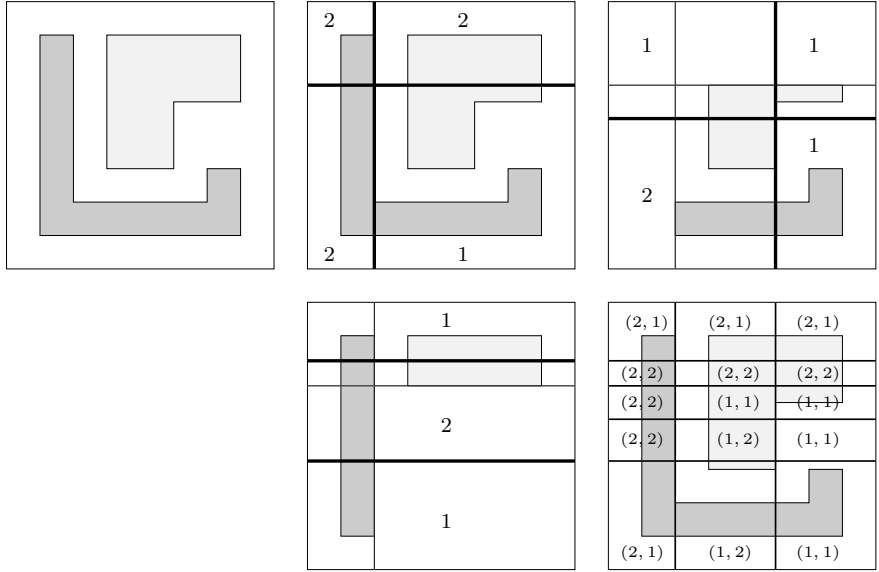
\begin{figure}[htbp]
  \centering
  \input{induction2.pstex_t}
  \caption{The reduction in action}
  \label{fig:induction}
\end{figure}\par
  First of all find numbers $m_i^j$,
  $j=0,\ldots , k_1$, such that
  \begin{itemize}
  \item $m_1^0+\cdots+m_d^0=n(k_1-1)$,
  \item $\sum_{i=1}^dm_i^j=n(k_2-1)$ for all $j=1,\ldots,k_1$, and
    \item $\sum_{j=0}^{k_1}m_i^j=m_i$ for all $j=1,\ldots,d$.
  \end{itemize}
This is certainly possible as
\begin{align*}
  n(k_1-1)+k_1n(k_2-1)=n(k-1)=\sum_{i=1}^dm_i.
\end{align*}
In the Figure we chose $m_1^0=m_2^0=1$, $m_1^1=m_2^1=1$, $m_1^2=0$,
and $m_2^2=2$.

By assumption there exists a set $\mathcal X^0$ of $n(k_1-1)$
hyperplanes of which $m_i^0$ are perpendicular to $e_i$ and
$\omega^0:C(\mathcal X^0)\rightarrow [k_1]$ such that
$\mu_i(A_j^0)=\frac{1}{k_1}$ for all $i=1,\ldots,n$, $j=1,\ldots,k_1$,
where $A_j^0$ are the unions of cells associated to $\omega^0$.

For $j=1,\dots,k_1$, consider the rescaled restrictions of the measures
to the regions $A_j^0$, i.e.,
\begin{align*}
  \mu_i^j(S)=k_1\mu_i(S\cap A_j^0).
\end{align*}
In other words for each $j$, $\mu_1^j,\ldots,\mu_n^j$ is a set
of $n$ probability measures on $I^d$ which have support only in
$A^0_j$.  Now for each $j$ let by assumption $\mathcal X^j$ be a set
of $n(k_2-1)$ hyperplanes of which $m_i^j$ are perpendicular to $e_i$
and $\omega^j:C(\mathcal X^j)\rightarrow[k_2]$ such that for the
associated $A_i^j$ we have
\begin{align*}
  \mu_i^j(A_{i'}^j)=\frac{1}{k_2}
\end{align*}
for all $j=1,\ldots,k_1$, and $i,i'=1,\ldots,k_2$.

We will now construct the desired pair $(\mathcal X,\omega)$ as
follows. Let $\mathcal X=X^0\cup X^1\cup\cdots\cup X^{k_1}$. To define
the map $\omega:C(\mathcal X)\rightarrow [k_1]\times[k_2]\cong [k]$ consider a
cell $c\in C(\mathcal X)$. Let $j_1\in[k_1]$ be the unique element
with $c\subseteq A^0_{j_1}$ and $j_2$ be the unique element in $[k_2]$
with $c\subseteq  A^{j_1}_{j_2}$. Then let
\begin{align*}
  \omega(c)=\left(\omega^0(A^0_{j_1}),\omega^{j_2}(A^{j_1}_{j_2})\right).
\end{align*}
\end{proof}
Applying the previous lemma we will now prove Theorem \ref{thm:main}.
\begin{proof}
  As we may assume $k$ to be prime, let $G=\mathbb Z_p$ be the cyclic
  group of prime order $p$. Let $\mathbb E$ be the space of all
  $n\times p$\,-matrices with row sums equal to zero.  $G$ acts on
  $\mathbb E$ by cyclic column permutations.  Let us construct a
  continuous $G$-equivariant map
  $f:\Omega(m_1,\ldots,m_d;G)\rightarrow \mathbb E$ such that each
  zero of this map corresponds to a desired solution. Let
  $(\mathfrak{t},\omega)=(t_1,\ldots,t_d,\omega)\in
  \Omega(m_1,\ldots,m_d;G)$. Following the notation from
  Section~\ref{sec:G-labelled polytopes}, each $t_i\in \Delta_{m_i}$
  is a partition $0=x_0^i\leq x_1^i\leq\ldots \leq x_{m_i}^i\leq
  x_{m_i+1}^i\leq 1$ of the unit interval $[0,1]$. Let
  $\mathcal{X}:=\cup_{i=1}^d~\mathcal{X}_i$ where $\mathcal{X}_i
  :=\{H^i_j\}_{j=1}^{m_i}$ is the collection of hyperplanes orthogonal
  to the unit vector $e_i$ defined by $H^i_j :=\{y\in \mathbb{R}^d\mid
  y_i = x^i_j\}$. The collection $\mathcal{X}$ dissects $I^d$ into
  $m_1\cdot\ldots\cdot m_d$ elementary $d$-parallelepipeds while the
  coloring function $\omega : C(\mathcal{X})\rightarrow [p]$ colors
  these $d$-parallelepipeds by $p$ colors. Let
  $A_i:=\cup~\omega^{-1}(i)$ be the union of all $d$-parallelepipeds
  colored by the color $i$. By construction, an element
  $(\mathfrak{t},\omega)$ corresponds to a fair division if
  $\mu_j(A_i)=1/p$ for each $i$ and $j$. Consequently the vector
  $v_i=v_i(\mathfrak{t},\omega):=(\mu_j(A_i)-1/p)_{j=1}^n\in
  \mathbb{R}^n$, which continuously depends on the input data
  $(\mathfrak{t},\omega)$, is equal to $0$ if and only if the division
  is fair from the point of view of $i$-th player (``thief''). By
  definition let
$$
f((\mathfrak{t},\omega)):= [v_1, v_2,\ldots,v_p]
$$
be the map $f: \Omega(\mathfrak{m};G)\rightarrow \mathbb{E}$
obtained by writing $v_i$ as column vectors of a matrix in
$\mathbb{E}$. The map $f$ is obviously $G$-equivariant. By
Theorem~\ref{thm:topo-shell} and
Corollary~\ref{cor:connectedness}, $\Omega =
\Omega(\mathfrak{m};G)$ is a $n(p-1)$-connected, free $G$-cell
complex. Hence, $\Omega, \mathbb{E}$ and $f$ together satisfy the
conditions of Theorem~\ref{thm:BSS-D}. Consequently $f$ must have
a zero which completes the proof of the theorem.
\end{proof}

\section{Concluding remarks}\label{sec:concluding}

It is customary to formulate consensus division theorems for
(vector-valued) measures $\mu$ that are continuous i.e.\ defined
by density functions $d\mu = f\cdot dm$, where $m$ is the Lebesgue
measure. It is not difficult to see that majority of these results
(including our Theorems~\ref{thm:2-dim-necklace} and
\ref{thm:main}) hold for much more general classes of measures.
For a broader perspective on this problem and other examples of
consensus division theorems the reader is referred to \cite{Ziv04,
MVZ}. Here we restrict ourselves to the observation that the
measures used in multidimensional splitting necklace theorems do
not have to be positive. Moreover, the continuity condition can be
replaced by a much weaker condition that $\mu(\partial(Q))=0$
where $Q\subset I^d$ is an arbitrary parallelepiped and
$\partial(Q)$ its boundary.

\medskip The ``rainbow complexes'' $\Omega(Q;[k])$, introduced in
Section~\ref{sec:G-labelled polytopes}, appear to have some
independent interest as topological/geometric objects which
capture some of the  combinatorial properties of the underlying
polytope $Q$. For example if $Q\subset \mathbb{R}^d$ is a
simplicial polytope, then the Euler characteristic
$\chi(\Omega(Q;[k]))$ is given by the formula
$$
\chi(\Omega(Q;[k])) = k\cdot F_Q(-k):= f_0k-f_1k^2+\ldots
+(-1)^{d-1}f_{d-1}k^d +(-1)^d k^{d+1},
$$
where $(f_0,f_1,\ldots,f_{d-1},f_d)$ is the $f$-vector of $Q$. A
broader outlook should place rainbow complexes $\Omega(Q;[k])$ and
their generalizations into the category of combinatorially defined
configuration spaces associated to polytopes, (Eulerian) posets,
simplicial complexes etc. In this generality they could be seen as
relatives of toric varieties and their combinatorial counterparts
(extensions) such as moment-angle complexes $\mathcal{Z}_K$
\cite{BP00, DJ91}, homotopy colimits over posets \cite{Crelle}
etc.

\medskip\noindent
{\bf Acknowledgement:} The authors would like to thank the
organizers of the special program ``Computational Applications of
Algebraic Topology'', hosted by the Mathematical Sciences Research
Institute (MSRI, Berkeley, Fall 2006), for the support, excellent
working conditions and stimulating research atmosphere. The second
author acknowledges the support by the Serbian Ministry of Science
(projects 144014 and 144026).

\end{document}

%% file: omegacell.pstex_t
\begin{picture}(0,0)%
\epsfig{file=omegacell.pstex}%
\end{picture}%
\setlength{\unitlength}{1381sp}%
\begingroup\makeatletter\ifx\SetFigFont\undefined%
\gdef\SetFigFont#1#2#3#4#5{%
  \reset@font\fontsize{#1}{#2pt}%
  \fontfamily{#3}\fontseries{#4}\fontshape{#5}%
  \selectfont}%
\fi\endgroup%
\begin{picture}(11925,5272)(9688,-8854)
\end{picture}%

%% file: induction2.pstex_t
\begin{picture}(0,0)%
\epsfig{file=induction2.pstex}%
\end{picture}%
\setlength{\unitlength}{1381sp}%
\begingroup\makeatletter\ifx\SetFigFont\undefined%
\gdef\SetFigFont#1#2#3#4#5{%
  \reset@font\fontsize{#1}{#2pt}%
  \fontfamily{#3}\fontseries{#4}\fontshape{#5}%
  \selectfont}%
\fi\endgroup%
\begin{picture}(15666,10287)(2389,-11194)
\put(10501,-1411){\makebox(0,0)[lb]{\smash{{\SetFigFont{8}{9.6}{\rmdefault}{\mddefault}{\updefault}{\color[rgb]{0,0,0}$2$}%
}}}}
\put(10201,-10411){\makebox(0,0)[lb]{\smash{{\SetFigFont{8}{9.6}{\rmdefault}{\mddefault}{\updefault}{\color[rgb]{0,0,0}$1$}%
}}}}
\put(10201,-8686){\makebox(0,0)[lb]{\smash{{\SetFigFont{8}{9.6}{\rmdefault}{\mddefault}{\updefault}{\color[rgb]{0,0,0}$2$}%
}}}}
\put(10201,-6811){\makebox(0,0)[lb]{\smash{{\SetFigFont{8}{9.6}{\rmdefault}{\mddefault}{\updefault}{\color[rgb]{0,0,0}$1$}%
}}}}
\put(13801,-1861){\makebox(0,0)[lb]{\smash{{\SetFigFont{8}{9.6}{\rmdefault}{\mddefault}{\updefault}{\color[rgb]{0,0,0}$1$}%
}}}}
\put(16801,-1861){\makebox(0,0)[lb]{\smash{{\SetFigFont{8}{9.6}{\rmdefault}{\mddefault}{\updefault}{\color[rgb]{0,0,0}$1$}%
}}}}
\put(16801,-3661){\makebox(0,0)[lb]{\smash{{\SetFigFont{8}{9.6}{\rmdefault}{\mddefault}{\updefault}{\color[rgb]{0,0,0}$1$}%
}}}}
\put(13801,-4561){\makebox(0,0)[lb]{\smash{{\SetFigFont{8}{9.6}{\rmdefault}{\mddefault}{\updefault}{\color[rgb]{0,0,0}$2$}%
}}}}
\put(8101,-5611){\makebox(0,0)[lb]{\smash{{\SetFigFont{8}{9.6}{\rmdefault}{\mddefault}{\updefault}{\color[rgb]{0,0,0}$2$}%
}}}}
\put(8101,-1411){\makebox(0,0)[lb]{\smash{{\SetFigFont{8}{9.6}{\rmdefault}{\mddefault}{\updefault}{\color[rgb]{0,0,0}$2$}%
}}}}
\put(10501,-5611){\makebox(0,0)[lb]{\smash{{\SetFigFont{8}{9.6}{\rmdefault}{\mddefault}{\updefault}{\color[rgb]{0,0,0}$1$}%
}}}}
\put(15001,-6811){\makebox(0,0)[lb]{\smash{{\SetFigFont{6}{7.2}{\rmdefault}{\mddefault}{\updefault}{\color[rgb]{0.000,0.000,0.000}$(2,1)$}%
}}}}
\put(16576,-8911){\makebox(0,0)[lb]{\smash{{\SetFigFont{6}{7.2}{\rmdefault}{\mddefault}{\updefault}{\color[rgb]{0.000,0.000,0.000}$(1,1)$}%
}}}}
\put(16576,-10936){\makebox(0,0)[lb]{\smash{{\SetFigFont{6}{7.2}{\rmdefault}{\mddefault}{\updefault}{\color[rgb]{0.000,0.000,0.000}$(1,1)$}%
}}}}
\put(13426,-8911){\makebox(0,0)[lb]{\smash{{\SetFigFont{6}{7.2}{\rmdefault}{\mddefault}{\updefault}{\color[rgb]{0.000,0.000,0.000}$(2,2)$}%
}}}}
\put(13426,-10936){\makebox(0,0)[lb]{\smash{{\SetFigFont{6}{7.2}{\rmdefault}{\mddefault}{\updefault}{\color[rgb]{0.000,0.000,0.000}$(2,1)$}%
}}}}
\put(13426,-8236){\makebox(0,0)[lb]{\smash{{\SetFigFont{6}{7.2}{\rmdefault}{\mddefault}{\updefault}{\color[rgb]{0.000,0.000,0.000}$(2,2)$}%
}}}}
\put(13426,-7711){\makebox(0,0)[lb]{\smash{{\SetFigFont{6}{7.2}{\rmdefault}{\mddefault}{\updefault}{\color[rgb]{0.000,0.000,0.000}$(2,2)$}%
}}}}
\put(16576,-8236){\makebox(0,0)[lb]{\smash{{\SetFigFont{6}{7.2}{\rmdefault}{\mddefault}{\updefault}{\color[rgb]{0.000,0.000,0.000}$(1,1)$}%
}}}}
\put(16576,-7711){\makebox(0,0)[lb]{\smash{{\SetFigFont{6}{7.2}{\rmdefault}{\mddefault}{\updefault}{\color[rgb]{0.000,0.000,0.000}$(2,2)$}%
}}}}
\put(15151,-8911){\makebox(0,0)[lb]{\smash{{\SetFigFont{6}{7.2}{\rmdefault}{\mddefault}{\updefault}{\color[rgb]{0.000,0.000,0.000}$(1,2)$}%
}}}}
\put(15151,-8236){\makebox(0,0)[lb]{\smash{{\SetFigFont{6}{7.2}{\rmdefault}{\mddefault}{\updefault}{\color[rgb]{0.000,0.000,0.000}$(1,1)$}%
}}}}
\put(15151,-7711){\makebox(0,0)[lb]{\smash{{\SetFigFont{6}{7.2}{\rmdefault}{\mddefault}{\updefault}{\color[rgb]{0.000,0.000,0.000}$(2,2)$}%
}}}}
\put(15001,-10936){\makebox(0,0)[lb]{\smash{{\SetFigFont{6}{7.2}{\rmdefault}{\mddefault}{\updefault}{\color[rgb]{0.000,0.000,0.000}$(1,2)$}%
}}}}
\put(13501,-6811){\makebox(0,0)[lb]{\smash{{\SetFigFont{6}{7.2}{\rmdefault}{\mddefault}{\updefault}{\color[rgb]{0.000,0.000,0.000}$(2,1)$}%
}}}}
\put(16576,-6811){\makebox(0,0)[lb]{\smash{{\SetFigFont{6}{7.2}{\rmdefault}{\mddefault}{\updefault}{\color[rgb]{0.000,0.000,0.000}$(2,1)$}%
}}}}
\end{picture}%